\documentclass[12pt,a4]{article}
\usepackage{amsmath}
\usepackage{amssymb}
\usepackage{amsfonts,amssymb}
\usepackage[]{graphicx}
\usepackage[cp1251]{inputenc}
\usepackage[english, russian]{babel}
\usepackage[a4paper, left=40mm, right=10mm, bottom=30mm, head=0mm]{geometry}

\begin{document}

\centerline{\textbf{\Large On the connections of an unbounded }}
\centerline{\textbf{\Large  Hille-Phillips   functional calculus with }}
\centerline{\textbf{\Large Bochner-Phillips  functional calculus}}

\

\centerline{A. R. Mirotin}

\centerline{amirotin@yandex.ru}

\

The extension of  Hille-Phillips functional calculus of semigroup generators  which leads to unbounded operators is considered. Connections of this calculus to Bochner-Phillips functional calculus are indicated. In particular, the multiplication rule and the composition rule are proved. Several examples are given.

\

Key wards: Hille-Phillips functional calculus, Bochner-Phillips functional calculus, fractional powers of operators, subordination.

\

\centerline{\textbf{О  СВЯЗЯХ НЕОГРАНИЧЕННОГО ФУНКЦИОНАЛЬНОГО   }}
\centerline{\textbf{ ИСЧИСЛЕНИЯ ХИЛЛЕ-ФИЛЛИПСА }}
\centerline{\textbf{С ИСЧИСЛЕНИЕМ БОХНЕРА-ФИЛЛИПСА}}

\

\centerline{\textbf{А. Р. Миротин}}

\

Дается расширение функционального исчисления Хилле-Филлипса генераторов $C_0$-полугрупп, изложенного в их известной монографии, приводящее к неограниченным операторам. Указаны связи этого исчисления с исчислением Бохнера-Филлипса. В частности, доказано правило произведения и теорема о сложной функции. Рассмотрены примеры.

\

Ключевые слова: функциональное исчисление Хилле-Филлипса, функциональное исчисление Бохнера-Филлипса, дробные степени операторов, подчиненная полугруппа.

\

\textbf{1. Введение.} В  монографии \cite{1}, гл.  XV -- XVI  построено функциональное исчисление генераторов $C_0$-полугрупп, использующее класс $L\mathcal{M}$ функций, представимых в виде преобразований Лапласа
\[
La(s):=\int _{0}^{\infty }e^{st} da(t) \; (s<0)
\]
$\sigma$-конечных комплексных регулярных борелевских мер $a$ на $\mathbb{R}_+$. Пространство таких мер будет обозначаться $\mathcal{M}(\mathbb{R}_+)$. При этом условия, налагаемые в \cite{1} на меру и полугруппу, приводят к тому, что возникающие в результате операторы ограничены. Там же (с. 463) поставлена задача построения расширения этого исчисления, приводящего к неограниченным операторам. В работе развивается подход к такому расширению, анонсированный ранее в  \cite{Izv}. Случай генераторов групп рассматривался в \cite{BHK}.
Подход, основанный на другом классе символов (причем для наборов нескольких генераторов), появился в \cite{LS}. Необходимость исчисления подобного типа вызвана также потребностями функционального исчисления Бохнера-Филлипса, использующего класс $\mathcal{T}$ отрицательных функций Бернштейна  (см., например, \cite{2} -- \cite{9}, а  также \cite{10}). Ниже будут установлены связи между этими исчислениями, в частности, доказаны правило умножения и теорема о сложной функции.
 Всюду ниже $A$ есть генератор ограниченной $C_0$-полугруппы $T$ в банаховом пространстве $X$  с областью определения  $D(A)$ и образом  $\mathrm{Im} A.$  Через   $LB(Y,X)$ обозначается  пространство линейных ограниченных операторов, действующих между банаховыми пространствами $Y$ и $X,$ $LB(X):=LB(X,X).$ Если $f$ есть функция на $\mathbb{R}_+,$ то через $\widehat{f}$  будет обозначаться преобразование Лапласа меры $f(t)dt.$ Для регулярной борелевской меры $\mu$ на $\mathbb{R}_+$ через $\mu(t)$ обозначается ее функция распределения, нормированная условиями $\mu(0)=0,$ $\mu(t)=(\mu(t-0)+\mu(t+0))/2$ при $t>0.$ Конец доказательства или примера обозначается знаком $\Box$.

\

\textbf{2. Основное определение.} Следующее определение  формально совпадает с определением, предложенным Хилле и Филлипсом в монографии \cite{1}, но мы отказываемся от  наложенных там ограничений, гарантирующих существование интеграла и ограниченность определяемого им оператора (см. также работу \cite{BHK}, посвященную генераторам групп).

\textbf{Определение 1. } Для функции $g$ из $L\mathcal{M},$ $g=La,$ $a\in \mathcal{M}(\mathbb{R}_+)$) положим

\[
g(A)x=\int\limits _{0}^{\infty }T(t)xda(t),
\]
где область определения $D_{0}(g(A))$ этого оператора состоит из тех $x\in X,$ для которых интеграл в правой части существует в смысле Бохнера.

\

Следующий пример иллюстрирует определение 1.

\textbf{Пример 1.} Пусть $g(s)=s^{-1}, s<0.$ Тогда  $g=La,$  где  $a=-\mathrm{mes}$ ($\mathrm{mes}$ --- мера Лебега на $\mathbb{R}_+$).  Таким образом, в силу определения 1 при $x\in D_{0} (g(A))$
$$
g(A)x=-\int\limits _{0}^{\infty }T(t)xdt.
$$
Предположим, что генератор $A$ инъективен, а полугруппа $T$ сильно устойчива (т.~е. $T(n)y\to 0$  при  $y\in X,$ $n\to \infty$) и покажем, что из определения 1 следует равенство $g(A)=A^{-1}.$ Пусть $x\in \mathrm{Im} A$ и $y=A^{-1}x.$ Тогда
$$
g(A)x=-\int\limits _{0}^{\infty }T(t)Aydt=-\int\limits _{0}^{\infty }dT(t)y=y.
$$
Таким образом, $D_{0} (g(A))\supseteq \mathrm{Im} A$ и при  $x\in \mathrm{Im} A$ имеем $g(A)x=A^{-1}x.$ Нам осталось доказать включение $D_{0} (g(A))\subseteq \mathrm{Im} A.$ С этой целью выберем произвольно $x\in D_{0} (g(A))$ и рассмотрим последовательность
$$
y_n:=-\int\limits _{0}^{n}T(t)xdt.
$$
Положим $y:=\lim\limits_{n\to\infty}y_n.$ Как известно, $y_n\in D(A)$  и
$$
Ay_n=-A\int\limits _{0}^{n}T(t)xdt=x-T(n)x.
$$
Поэтому $Ay_n \to x$ $(n\to\infty),$ и  в силу замкнутости оператора $A$ имеем $y\in D(A)$ и  $Ay=x,$ что и завершает доказательство.$\Box$

\

\textbf{3. Теоремы о замкнутости $g(A)$.} Прежде всего, нас интересуют условия, при которых оператор $g(A)$ будет плотно определен и замкнут (замыкаем).

\

 \textbf{Лемма 1. }1) \textit{Если $\mathrm{Im} A\subset D_{0} (g(A))$, то оператор $g(A)A$ ограничен относительно $A$;}

 2) \textit{если дополнительно предположить, что оператор $A$ инъективен, то $g(A)$ замкнут на подпространстве $D(A)\cap D_{0} (g(A))$, наделенном нормой графика;}

 3) \textit{если $\int _{0}^{\infty }\left\|T(t)\right\|d\left|a\right|(t)<\infty$, то оператор $g(A)$ ограничен на $X,$ и рассматриваемое исчисление согласовано с классическим исчислением Хилле-Филлипса.}

\

 Доказательство. 1) Для любого $x\in D(A)$ определим операторы

\[
B_{n} x=\int _{0}^{n}T(t)Axda(t) ,
\]
(интеграл существует в смысле Бохнера, так как функция $r\mapsto \left\| T(t)\right\| $ ограничена на [0,\textit{n}] по принципу равномерной ограниченности). Поскольку у нас $Ax\in D_{0} (g(A))$, то $B_{n} x\to g(A)Ax\; (n\to \infty )$ при всех $x\in D(A)$. Далее, так как оператор $A$ замкнут, то пространство $Y=D(A)$, наделенное нормой графика $\left\| x\right\| _{Y} =\left\| x\right\| +\left\| Ax\right\| $, банахово. Кроме того, $B_{n} \in LB(Y,X)$, поскольку

\[
\left\| B_{n} x\right\| \le \left(\int _{0}^{n}\left\|T(t)\right\|  d|a|(t)\right)\left\| x\right\| _{Y} .
\]

В силу теоремы Банаха-Штейнгауза оператор $g(A)A$ тоже принадлежит $LB(Y,X)$, а потому ограничен относительно $A.$

 2) Заметим сначала, что при $x\in D(A)\cap D_{0}(g(A))$ справедливо равенство
 $$
 Ag(A)x=g(A)Ax.  \eqno(1)
$$
Действительно, с учетом замкнутости $A$ и сходимости интегралов имеем

\[
g(A)Ax=\int _{0}^{\infty }T(t)Axda(t) =\int _{0}^{\infty }AT(t)xda(r) =Ag(A)x,
\]
поскольку $\mathrm{Im}A\subset D_{0} (g(A))$. Теперь, если $A$ инъективен, то с помощью утверждения 1) получаем, что оператор $g(A)x=A^{-1} g(A)Ax$ замнут на подпространстве $D(A)\cap D_{0} (g(A))$ пространства $Y$ как произведение замкнутого и ограниченного операторов.

 3) Это следует из свойств интеграла Бохнера.  $\Box$

\

Всюду далее для функции $f$ на $\mathbb{R}_+$ через $\widehat{f}$ обозначается ее преобразование Лапласа, т.~ е.
\[
\widehat{f}(s):=\int _{0}^{\infty }e^{st}f(t)dt \; (s<0).
\]

\

\textbf{ Следствие 1.} \textit{Пусть $g=\widehat{f}$ есть преобразование Лапласа функции $f,$
причем при $x$ из $D(A)$ существуют $\mathop{\lim }\limits_{n\to \infty } f(n)T(n)x$ и $\int _{0}^{\infty }T(t)xdf(t) $. Тогда $\mathrm{Im}A\subset D_{0} (g(A))$ и справедливы все утверждения леммы 1.}

\

Доказательство. Интегрируя по частям, получаем при всех $x\in D(A)$
$$
\int _{0}^{n}T(t)Axf(t)dt= \int _{0}^{n}f(t)dT(t)x=f(n)T(n)x-f(0)x-\int _{0}^{n}T(t)xdf(t),  \eqno(2)
$$
причем правая часть имеет предел при $n\to \infty $.    $\Box$

\

\textbf{ Следствие 2. }\textit{В условиях части 1 леммы 1 оператор $g(A)A$ ограничен вместе с $A.$
}

\

\textbf{ Следствие 3. }\textit{Если в условиях части 1 леммы 1 оператор \textit{А} инъективен, то $g(A)|\mathrm{Im}A$ ограничен относительно $A^{-1}$.}

\

\textbf{ Следствие 4. }\textit{Если в условиях части 1 леммы 1 существует ограниченный обратный оператор $A^{-1}$ на $\mathrm{Im}A$, то $g(A)$ ограничен на $\mathrm{Im}A$.
}

\

\textbf{ Теорема 1. }\textit{Пусть $g=\widehat{f}$, где функция $f$ такова, что оператор
$$
Sx:=\int\limits_{0}^{\infty }T(t)xdf(t)
 $$
  ограничен, а последовательность $f(n)T(n)$  сходится на $D(A)$ сильно к оператору $B\in LB(D(A)).$ Тогда}

 1) $\mathrm{Im}A\subset D_{0} (g(A))$\textit{ и оператор $g(A)A$ ограничен на }$D(A)$;

 2) \textit{если оператор $A$  инъективен, то оператор $g(A)|D(A)\cap D_{0}(g(A))$ замкнут, а если еще $B=0,$ то замкнут также и оператор} $g(A)|\mathrm{Im}(A)$.

\

 Доказательство. 1) Включение  $\mathrm{Im}A\subset D_{0}(g(A))$ сразу вытекает из следствия 1. Переходя к пределу в формуле (2), получаем при $x\in D(A)$

$$
g(A)Ax=\int _{0}^{\infty }T(t)Axf(t)dt=B x-f(0)x-\int _{0}^{\infty }T(t)xdf(t), \eqno(3)
$$
откуда и следует ограниченность $g(A)A$ на $D(A)$.

 2) Здесь первое утверждение следует из 1) и равенства $g(A)x=A^{-1} g(A)Ax$ ($x\in D(A)\cap D_{0}(g(A))$) как в доказательстве леммы 1. Пусть теперь  $B=0$. Полагая в (3) $y=Ax,$ имеем
\[
g(A)y=-f(0)A^{-1} y-\int _{0}^{\infty }T(t)A^{-1} ydf(t) .
\]
Но оператор
\[
y\mapsto A\int _{0}^{\infty }T(t)A^{-1} ydf(t) =\int _{0}^{\infty }AT(t)A^{-1} ydf(t) =\int _{0}^{\infty }T(t)ydf(t)=Sy
\]
ограничен на $\mathrm{Im}A,$ а потому оператор
\[
g(A)y=-A^{-1} \left(f(0)y+Sy\right)
\]
замкнут на  $\mathrm{Im} A$ как произведение замкнутого и ограниченного операторов.   $\Box$

\

Отрицательные дробные степени неограниченных операторов  рассматривались многими авторами (см., например, \cite{Kr}, \cite{Hen}, \cite{MS}, \cite{KRek}). Определение 1 приводит к определению отрицательных дробных степеней генераторов полугрупп, формально совпадающему с определением из \cite[c. 32--33]{Hen}. \footnote{В \cite{Hen} отрицательные дробные степени определяются при условиях, что   оператор $-A$  секториальный и $\mathrm{Re}\sigma(A)<0,$ и являются ограниченными операторами.}

\textbf{Пример 2}. Пусть $g(s)=(-s)^{-\alpha}, s<0, \alpha>0.$ Тогда  $g=\widehat{f},$  где  $f(r)=1/\Gamma(\alpha)r^{\alpha-1}.$   Следовательно,  мы можем в соответствии с определением 1 положить
$$
(-A)^{-\alpha}x:=\frac{1}{\Gamma(\alpha)}\int\limits_{0}^{\infty}T(t)xt^{\alpha-1}dt,
$$
считая, что $D_0((-A)^{-\alpha})$ состоит из тех $x\in X,$ при которых интеграл существует в смысле Бохнера. Предположим, что  $\alpha>1,$ оператор $A$ инъективен, а $C_0$-полугруппа $T$ удовлетворяет оценке $\|T(t)\|\le C/t^\delta$ с  константами  $\delta>\alpha-1,$ $C>0.$
Тогда выполнены все условия  теоремы 1, а потому  $\mathrm{Im}A\subset D_0((-A)^{-\alpha})$ и оператор $(-A)^{-\alpha}|\mathrm{Im} A$ замкнут.

Если же $0<\alpha<1,$ то мы можем при $\delta>\alpha$ положить $(-A)^{-\alpha}:=(-A)^{-(1+\alpha)}(-A).$ В этом случае $D((-A)^{-\alpha})=D(A),$ и при $x\in D(A),$ интегрируя по частям, для $(-A)^{-\alpha}x$ получим ту же формулу, что и выше. В самом деле, тогда $t^\alpha T(t)x\to 0$ ($t\to\infty$), а потому
$$
(-A)^{-\alpha}x:=\frac{1}{\Gamma(\alpha+1)}\int\limits_{0}^{\infty}T(t)(-Ax)t^{\alpha}dt=
$$
$$
-\frac{1}{\alpha\Gamma(\alpha)}\left(t^\alpha T(t)x\left|_{0}^\infty\right.-\alpha\int\limits_{0}^{\infty}T(t)xt^{\alpha-1}dt\right)=
\frac{1}{\Gamma(\alpha)}\int\limits_{0}^{\infty}T(t)xt^{\alpha-1}dt.
$$

  Если дополнительно предположить, что полугруппа   сжимающая, то $(-A)^{-\alpha}x\to x $ при $x\in D(A),$ $\alpha\to +0.$ Действительно, в этом случае $\|T(t)x-e^{-t}x\|\le t\|Ax-x\|$  (см., например, \cite[гл. 1, лемма 2.9]{Goldst}). Так как
  $$
  (-A)^{-\alpha}x-x=\frac{1}{\Gamma(\alpha)}\int\limits_{0}^{\infty}(T(t)x-e^{-t}x)t^{\alpha-1}dt,
  $$
  то
  $$
  \|(-A)^{-\alpha}x-x\|\le\frac{1}{\Gamma(\alpha)}\left(\int\limits_{0}^{1}\|T(t)x-e^{-t}x\|t^{\alpha-1}dt+
 \int\limits_{1}^{\infty}\|T(t)x-e^{-t}x\|t^{\alpha-1}dt \right)\le
  $$
  $$
  \frac{1}{\Gamma(\alpha)}\left(\int\limits_{0}^{1}\|Ax-x\|t^{\alpha}dt+
 \int\limits_{1}^{\infty}\left(\frac{M}{t^\delta}+e^{-t}\right)\|x\|t^{\alpha-1}dt \right)\le
  $$
  $$
  \frac{1}{\Gamma(\alpha)}\left(\frac{\|Ax-x\|}{\alpha+1}+\frac{M\|x\|}{\delta-\alpha}+\|x\|e^{-1}\right)\to 0
  $$
при $\alpha\to +0. \Box$

\

\textbf{4. Связь с исчислением Бохнера-Филлипса.} Следующие теоремы устанавливают связь между рассматриваемым исчислением и исчислением Бохнера-Филлипса  (относительно последнего см., например, \cite{3}, \cite{5}). Ниже через  $\mathcal{T}$ будет обозначаться класс отрицательных функций Бернштейна одного переменного.  Функция $\psi$ из $\mathcal{T}$ аналитична в левой полуплоскости и допускает интегральное представление
 $$
 \psi (z)=c_0+\int _{0}^{\infty }\left(e^{zu} -1\right)u^{-1} d\rho(u) \; (\mathrm{Re}z<0), \eqno(4)                             $$
где $c_0=\psi (0),$ $\rho$ --- положительная мера на $\mathbb{R}_+$, причем  $\int _{0}^{r}d\rho(u)<\infty,$ $\int _{r}^{\infty }u^{-1}d\rho(u)<\infty $ при $r>0.$

\

\textbf{Определение 2.} Для неположительной функции Бернштейна $\psi$ с интегральным представлением (4) и генератора  $A$ ограниченной
$C_0$-полугруппы $T$ на банаховом пространстве $X$ ее значение на операторе  $A$ при $x\in D(A)$ определяется интегралом Бохнера
$$
\psi(A)x = c_0x + \int\limits_{\mathbb{R}_+} (T(u)-I)xu^{-1}d\rho(u).
$$
При этом замыкание этого  оператора,  также обозначаемое
$\psi(A)$,  есть генератор ограниченной
$C_0$-полугруппы $g_t(A)$ на $X.$ Здесь $g_t(s):= e^{t\psi(s)},$ а
 оператор $g_t(A)$ понимается в смысле исчисления Хилле-Филлипса, поскольку функция $g_t(s)$ абсолютно монотонна на $(-\infty,0],$ а потому есть преобразование Лапласа некоторой субвероятностной меры $\nu_t$ (полугруппа $g_t(A)$ называется подчиненной полугруппе  $T$).

\

 \textbf{Теорема 2. }\textit{Пусть    $\psi \in \mathcal{T},$ $\psi(0)=0$. Тогда функция $\tilde{\psi }(s):=\psi (s)/s$ принадлежит $L\mathcal{M}$ и при всех $x\in D(A)\cap D_{0} (\tilde{\psi }(A))$ справедливо равенство}
$$
 \psi (A)x=A\tilde{\psi }(A)x, \eqno(5)
$$
\textit{где $\tilde{\psi }(A)$ понимается в смысле определения 1, а $\psi (A)$ --- в смысле исчисления Бохнера-Филлипса.}

\

 Доказательство. В силу формулы (4) и теоремы Фубини

\[
\begin{array}{l}
 {\tilde{\psi }(s)=\int\limits _{0}^{\infty }\left(\frac{e^{su} -1}{s} \right)u^{-1} d\rho (u) =\int\limits _{0}^{\infty }\left(\int\limits _{0}^{\infty }1_{[0;u]} (r)e^{sr} dr \right) u^{-1} d\rho (u)=} \\ { \quad  \quad \quad \quad \;  \; \; \; \; \quad  \quad  \quad \int\limits _{0}^{\infty }e^{sr} \left(\int\limits _{0}^{\infty }1_{[0;u]} (r) u^{-1} d\rho (u)\right)dr=\widehat{f}(s) ,}
  \end{array}
    \]
где $f(r)=\int _{r}^{\infty }u^{-1} d\rho (u),\;  $ а $1_{A}$ --- индикатор множества $A.$

 Следовательно, если $x\in D(A)\cap D_{0} (\tilde{\psi }(A))$, то
$$
\tilde{\psi }(A)x=\int\limits _{0}^{\infty }T(t)xf(t)dt. \eqno(6)
$$

 С другой стороны, по теореме Фубини для интеграла Бохнера
\[
\begin{array}{l}
 {\psi (A)x:=\int\limits _{0}^{\infty }\left(T(u)-I\right)xu^{-1} d\rho(u) =\int\limits _{0}^{\infty }\left(\int\limits _{0}^{u}AT(t)xdt \right)u^{-1} d\rho(u) =} \\ {\quad \quad \quad \int\limits _{0}^{\infty }\left(\int\limits _{0}^{\infty }1_{[0;u]} (t)AT(t)xdt \right)u^{-1} d\rho (u) =A\int\limits _{0}^{\infty }T(t)x\left(\int\limits _{0}^{\infty }1_{[0;u]}(t) u^{-1} d\rho(u)\right)dt=} \\ {\quad \quad \quad A\tilde{\psi }(A)x}
  \end{array}
  \]
(теорема Фубини применима, поскольку интеграл  Бохнера в (6) сходится).          $\Box$

\

 \textbf{Следствие 5. }\textit{Если\textbf{ }оператор $A$ имеет ограниченный обратный, то оператор $\tilde{\psi }(A)$ замкнут на подпространстве $D(A)\cap D_{0} (\tilde{\psi }(A))$.}

\

 Доказательство. В силу теоремы 4.1 из \cite{4} и ее следствия формула (5) влечет равенство $\tilde{\psi }(A)x=\psi (A)A^{-1} x$ ($x\in D(A)\cap D_{0} (\tilde{\psi }(A))$), правая часть которого есть произведение замкнутого и ограниченного операторов.       $\Box$

\

 \textbf{Следствие 6. }\textit{Оператор $\tilde{\psi}(A)$ отображает $D(A)\cap D_{0}(\tilde{\psi }(A))$ в $D(A)$.}

\

\textbf{ Следствие 7.} \textit{Оператор $A$ отображает $D(A)\cap D_{0} (\tilde{\psi }(A))$ в $D_{0} (\tilde{\psi }(A))$.}

\

Для доказательства следующих двух теорем нам понадобится такой вариант   теоремы $X_{38}$ из \cite[Введение]{Wiener}.

\textbf{Лемма 2.} \textit{Пусть комплекснозначная функция $\phi$ непрерывна на $(0,\infty),$  функция $a(u)$  монотонно возрастает и имеет ограниченную вариацию на любом интервале $(c,d)\subset \mathbb{R}_+,$ а функция  $\alpha(u,v)$ непрерывна по $u$ при каждом  $v>0,$ монотонно возрастает по $v$ и на любом интервале $(a,b)\subset \mathbb{R}_+,$ имеет ограниченную вариацию по $v$ равномерно по переменной $u,$ пробегающей любой интервал $(c,d)\subset \mathbb{R}_+.$ Тогда, если один из интегралов}
$$
\int\limits_{0}^{\infty }\left(\int\limits_{0}^{\infty }\phi(v)d_v\alpha(u,v)\right)da(u), \int\limits_{0}^{\infty }\phi(v)d_v\left(\int\limits_{0}^{\infty }\alpha(u,v)da(u)\right)
$$
\textit{существует и конечен при замене  $\phi$ на  $|\phi|,$ то эти интегралы равны.}

\

Доказательство. Если функция $\phi$ неотрицательна, то, в силу упоминавшейся теоремы $X_{38},$
для любого натурального $n$
$$
\int\limits_{0}^{\infty }\left(\int\limits_{0}^{n}\phi(v)d_v\alpha(u,v)\right)da(u)= \int\limits_{0}^{n}\phi(v)d_v\left(\int\limits_{0}^{\infty }\alpha(u,v)da(u)\right).
$$
В этом случае для доказательства достаточно положить $n\to\infty$ и применить теорему Б. Леви.

Случай комплекснозначной функции $\phi$ сводится к предыдущему, если заметить, что $\phi$ можно представить в виде $\phi=(\phi_1-\phi_2)+i(\phi_3-\phi_4)$, где все функции $\phi_k$ удовлетворяют неравенствам  $0\le \phi_k\le|\phi|$ и непрерывны. $\Box$

\

 Положим

\[
L\mathcal{M}_{T} :=\left\{g\in L\mathcal{M}: g=La,\forall x\in D(A)\, \mathop{\lim }\limits_{t\to \infty }a(t)\|T(t)x\|=0\, \right\}.
\]

\

 \textbf{Теорема 3} (ср. \cite[теорема 1]{5}).  \textit{Пусть   $T$ --- ограниченная полугруппа,  $g=La,$ где $a$ --- непрерывная положительная мера, $g\in LM_{T} $, $\psi \in \mathcal{T}.$ Тогда $h:=g\psi\in L\mathcal{M}_{T},$
 $D_{0}(h(A))\supset D(A)\cap D_{0}(g(A))$ и при $x\in D(A)\cap D_{0}(g(A))$ справедливы равенства}
\[
h(A)x=\psi(A)g(A)x=g(A)\psi(A)x.
\]

\

 Доказательство. Пусть $g=La,$  $a(t)$ --- такая (непрерывная)  функция распределения меры $a,$ что $a(0)=0.$ Тогда $h=Lb,$ где $b$ --- мера на $\mathbb{R}_+$ с функцией распределения
$$
 b(t)=\psi (0)a(t)+\int\limits_{0}^{\infty }\left(a(t-u)-a(t)\right)u^{-1} d\rho(u).\eqno(7)
$$

В самом деле, в силу леммы 2
$$
Lb(s)=\psi(0)La(s)+\int\limits_{0}^{\infty }e^{st}d_t\left(\int\limits_{0}^{\infty }(a(t-u)-a(t))u^{-1} d\rho(u)\right)=
$$
$$
\psi(0)g(s)+\int\limits_{0}^{\infty }u^{-1} d\rho(u)\int\limits_{0}^{\infty }e^{st}d_t\left(a(t-u)-a(t)\right)=
$$
$$
 \psi(0)g(s)+\int\limits_{0}^{\infty }\left(\int\limits_{0}^{\infty }e^{st}d_ta(t-u)-\int\limits_{0}^{\infty }e^{st}da(t)\right) u^{-1} d\rho(u)=
 $$
 $$
 \psi(0)g(s)+\int\limits_{0}^{\infty }\int\limits_{0}^{\infty }\left(e^{s(t+u)}-e^{st}\right)da(t) u^{-1} d\rho(u)=
$$
$$
\psi(0)g(s)+\int\limits_{0}^{\infty }e^{st}da(t)\int\limits_{0}^{\infty }\left(e^{us}-1\right)u^{-1} d\rho(u)=g(s)\psi(s).
$$

Далее, поскольку при $x\in \in D(A), t>0$
$$
b(t)T(t)x=\psi(0)a(t)T(t)x+\int\limits_{0}^{\infty }\left(a(t-u)-a(t)\right)T(t)xu^{-1} d\rho(u)=
$$
$$
\psi(0)a(t)T(t)x+\int\limits_{0}^{\infty }a(t)(T(t+u)x-T(t)x)u^{-1} d\rho(u)=a(t)T(t)\psi(A)x,
$$
то $h\in L\mathcal{M}_{T}.$ Кроме того, полагая  $t\to  +0$, получаем $b(+0)=0.$

  С помощью интегрирования по частям легко проверить, что при $x\in D(A)\cap D_{0}(g(A))$
 $$
 g(A)x=\int _{0}^{\infty }T(t)(-Ax)a(t)dt.\eqno(8)
  $$
 Следовательно, при этих $x$ справедливы равенства
\[
\begin{array}{l}
 {(T(u)-I)g(A)x=\int\limits _{0}^{\infty }T(u+t)(-Ax)a(t)dt -\int\limits_{0}^{\infty }T(t)(-Ax)a(t)dt =} \\ {\quad \quad \quad \quad \quad \quad \quad =\int\limits_{0}^{\infty }T(t)(-Ax)(a(t-u)-a(t))dt.} \end{array}
\]
Поэтому и с учетом (7) при $x\in D(A)\cap D_{0}(g(A))$ имеем (не нарушая общности, можно считать $\psi(0)=0$)
\[
 \psi(A)g(A)x=
 \]
 \[
 \begin{array}{l}
{\int\limits_{0}^{\infty }\left(T(u)-I)\right) g(A)xu^{-1}d\rho(u)=
\int\limits_{0}^{\infty }\left(\int\limits_{0}^{\infty }T(t)(-Ax)(a(t-u)-a(t))dt \right) u^{-1} d\rho(u)=} \\ {\int\limits_{0}^{\infty }T(t)(-Ax)\left(\int\limits_{0}^{\infty }(a(t-u)-a(t)) u^{-1} d\rho(u)\right)dt= \int\limits_{0}^{\infty }T(t)(-Ax)b(t)dt=-\int\limits_{0}^{\infty }b(t)dT(t)x=} \end{array}
\]
\[
-b(t)T(t)x|_{t=0}^\infty+\int\limits_{0}^{\infty }T(t)xdb(t)=h(A)x.
\]
 Значит, $D_{0}(h(A))\supset D(A)\cap D_{0}(g(A))$ и первое равенство доказано.

 Наконец заметим, что операторы $T(u)$ и $g(A)$ коммутируют в том смысле, что при $x\in D_0(g(A))$
 $$
 T(u)g(A)x=\int\limits_{0}^{\infty }T(u+t)xda(t)=\int\limits_{0}^{\infty }T(t)T(u)xda(t)=g(A)T(u)x,
 $$
  а потому при $x\in D(A)\cap D_0(g(A))$
$$
\begin{array}{l} {h(A)x=\psi (A)g(A)x=\int\limits_{0}^{\infty }g(A)\left(T(u)-I\right)xu^{-1} d\rho (u)= } \\ {\quad \quad \quad \quad \quad \quad \quad \quad =g(A)\int\limits_{0}^{\infty }\left(T(u)-I\right)xu^{-1} d\rho (u)= g(A)\psi (A)x.} \\ {} \end{array}\Box
$$

\

\textbf{Следствие 8.}  \textit{ Пусть    $g\in LM_{T}.$ Тогда $h(s):=sg(s)\in L\mathcal{M}_{T},$
 $D_{0}(h(A))\subset D(A)\cap D_{0}(g(A))$ и при $x\in D(A)\cap D_{0}(g(A))$ справедливы равенства}
\[
h(A)x=Ag(A)x=g(A)Ax.
\]

В самом деле, если в (4) в качестве $\rho$ взять меру Дирака, то при $c_0=0$ получим $\psi(s)=s.$

\

\textbf{Пример 3.} Пусть $g(s)=(-s)^{-\alpha}, s<0,$  $0<\alpha<1,$ а $C_0$-полугруппа $T$ удовлетворяет оценке $\|T(t)\|\le C/t^\delta$ с  константами  $\delta>\alpha,$ $C>0$ как в примере 2. И пусть $\psi(s)=-(-s)^{\beta}, s<0,$ $0<\beta<\alpha.$ Как известно, $\psi\in \mathcal{T}$ и
$$
-(-s)^{\beta}=\frac{\beta}{\Gamma(1-\beta)}\int\limits_{0}^{\infty }(e^{st}-1)t^{-\beta-1}dt.
$$
Тогда выполнены все условия теоремы 3, причем в силу примера 2 $D_0(g(A))=D_0(h(A))=D(A),$ а потому для генератора $A$ полугруппы $T$ при $x\in D(A)$ справедливо равенство
$$
(-A)^\beta(-A)^{-\alpha}x=(-A)^{\beta-\alpha}x. \Box
$$

\

В связи с идущим ниже следствием 8 отметим, что, если $\psi \in \mathcal{T}$, то функция $-1/\psi$ на $(-\infty,0)$ абсолютно монотонна, а потому принадлежит $L\mathcal{M}.$

\

 \textbf{Следствие 9.} \textit{Пусть  $\psi \in \mathcal{T}$, функция $1/\psi $ принадлежит $L\mathcal{M}_{T} $, и оператор
 $(1/\psi)(A)$ (в смысле определения 1) определен и ограничен на $X.$ Тогда оператор $\psi(A)$ обратим, и $\psi(A)^{-1} =
 (1/\psi)(A)$.}

Доказательство. Так как  $(1/\psi)(s)\psi(s)=1,$ то  силу теоремы 3 при $x\in D(A)$  имеем
$$
(1/\psi)(A)\psi(A)x=\psi(A)(1/\psi)(A)x=x.
$$
Поскольку оператор $\psi(A)(1/\psi)(A)$ замкнут как произведение замкнутого и ограниченного операторов, последнее равенство верно при всех  $x\in X.$ Следовательно, оператор  $\psi(A)$ биективен и  $\psi(A)^{-1} =
(1/\psi)(A).$
$\Box$

\

\textbf{Пример 4.} Пусть $\psi(s)=-\log(1-s).$ Известно, что $\psi\in \mathcal{T}$ (см., напр., \cite{3}). Кроме того, известно, что для функции
$$
\nu(t,-1)=\int\limits_{1}^{\infty }\frac{t^{\xi}}{\xi\Gamma(\xi)}d\xi
$$
справедливо равенство  $1/\log(-x)=\widehat{\nu(t,-1)}(x)$ при $x<0$
\cite[глава V, \S 5.7, (11)]{BE}. Следовательно, $1/\psi=\widehat{(-f)},$ где $f(t)=e^{-t}\nu(t,-1).$ Пусть полугруппа $T$ равномерно устойчива, т. е. $\|T(t)\|\le Me^{\omega t},$ где
$\omega<0.$ С помощью правила Лопиталя легко проверить, что $1/\psi\in L\mathcal{M}_T.$ Кроме того, оператор
 $(1/\psi)(A)$  определен и ограничен на $X,$ так как
 $$
 \|(1/\psi)(A)x\|\le \int\limits_{0}^{\infty }\|T(t)x\|f(t)dt\le M\int\limits_{0}^{\infty }e^{\omega t}f(t)dt\|x\| =\frac{M}{\log(1-\omega)}\|x\|.
 $$
 Таким образом, по следствию 8 существует ограниченный обратный оператор
 $$
 (\log(I-A))^{-1}x= \int\limits_{0}^{\infty }T(t)xe^{-t}\nu(t,-1)dt  \ (x\in X).\Box
 $$

\

\textbf{ Замечание 1.} В \cite{Trudy} для случая, когда $A$  есть негативный оператор (в смысле Коматсу)  в $X,$
доказана формула
 $$
 (\log(I-A))^{-1}x= \int\limits_{1}^{\infty }\frac{R(t,A)xdt}{\pi^2+\log^2(t-1)} \ (x\in X).
 $$

\

Положим
$$
X_T:=\{x\in X: \exists M(x)>0, \omega(x)<0\ \forall t\in \mathbb{R}_+\  \|T(t)x\|\le M(x)e^{\omega(x) t} \}.
$$

\

 \textbf{Теорема 4}.  \textit{Пусть     $h=La,$ где $a$ --- положительная мера,  $\psi \in \mathcal{T}.$ Тогда $h\circ\psi\in L\mathcal{M},$
  и при $x\in X_T$ справедливо равенство}
\[
(h\circ \psi)(A)x=h(\psi(A))x.
\]

\

 Доказательство. Справедлива формула
 $$
 (h\circ \psi)(s)=\int\limits_{0}^{\infty }e^{u\psi(s)}da(u).
 $$
Но, как  отмечалось выше,
$$
g_u(s)= e^{u\psi(s)}=\int\limits_{0}^{\infty }e^{sv}d\nu_u(v), u\ge 0
$$
для некоторой субвероятностной положительной меры $\nu_u$ на $\mathbb{R}_+,$ а потому
$$
 (h\circ \psi)(s)=\int\limits_{0}^{\infty }\left(\int\limits_{0}^{\infty }e^{sv}d\nu_u(v)\right)da(u).\eqno(9)
$$
Чтобы применить лемму 2, нам достаточно проверить, что функция $u\mapsto \nu_u(v)$  непрерывна на  $(0,\infty).$ Для этого воспользуемся 
теоремой непрерывности для преобразования Лапласа (см., например, \cite[глава XIII, \S 1, теорема 2а]{Fel}). Пусть $u_n, u>0$ и $u_n\to u \ (n\to\infty).$ Поскольку $L\nu_{u_n}(-\lambda)=e^{u_n\psi(-\lambda)}\to L\nu_u(-z)$ при $\lambda>0,$ то по указанной теореме $\nu_{u_n}\{J\}\to \nu_{u}\{J\}$ для любого конечного интервала  $J$ непрерывности функции $\nu_u,$ что позволяет сделать вывод о непрерывности функции $u\mapsto \nu_u(v)$ на  $(0,\infty).$

Стало быть, в силу леммы 2, из (9) следует, что
$$
 (h\circ \psi)(s)=\int\limits_{0}^{\infty }e^{sv}d_v\left(\int\limits_{0}^{\infty }d\nu_u(v)da(u)\right).
$$
Таким образом, $h\circ \psi\in L\mathcal{M}$ и при $x\in D_0((h\circ \psi)(A))$
$$
 (h\circ \psi)(A)x:=\int\limits_{0}^{\infty }T(v)xd_v\left(\int\limits_{0}^{\infty }d\nu_u(v)da(u)\right).\eqno(10)
$$
Заметим, что при $x\in X_T$ интеграл, стоящий в правой части формулы  (10), существует в смысле Бохнера, поскольку в силу леммы 2 и формулы (9)
$$
\int\limits_{0}^{\infty }\|T(v)x\|d_v\left(\int\limits_{0}^{\infty }d\nu_u(v)da(u)\right)\le
M(x)\int\limits_{0}^{\infty }e^{\omega(x)v}d_v\left(\int\limits_{0}^{\infty }d\nu_u(v)da(u)\right)=
$$
$$
M(x)\int\limits_{0}^{\infty }\left(\int\limits_{0}^{\infty }e^{\omega(x)v}d_v\nu_u(v)\right)da(u)=M(x)(h\circ \psi)(\omega(x)).
$$
Значит, $X_T\subset D_0((h\circ \psi)(A)).$

Далее, для любого функционала $\Lambda\in X'$  при $x\in X_T$  из (10) следует с учетом леммы 2, что
$$
 \Lambda((h\circ \psi)(A)x)=\int\limits_{0}^{\infty }\Lambda(T(v)x)d_v\left(\int\limits_{0}^{\infty }d\nu_u(v)da(u)\right)=
$$
$$
\int\limits_{0}^{\infty }\left(\int\limits_{0}^{\infty }\Lambda(T(v)x)d_v\nu_u(v)\right)da(u)=
\Lambda\left(\int\limits_{0}^{\infty }\int\limits_{0}^{\infty }T(v)x)d_v\nu_u(v)da(u)\right)=
$$
$$
\Lambda\left(\int\limits_{0}^{\infty }g_u(A)xda(u)\right)=\Lambda\left(h(\psi(A))x\right),
$$
откуда и следует доказываемое равенство. $\Box$

\

\textbf{Следствие 10.} \textit{Если полугруппа $T$ равномерно устойчива (т.~е. $\|T(t)\|\to 0, t\to\infty$), то $(h\circ \psi)(A)=h(\psi(A)),$ причем этот оператор ограничен.
}

\

\textbf{Пример 5.} Пусть $h(t)=(-t)^{-\alpha},$ $\psi(s)=-(-s)^\beta,$ $\alpha, \beta\in (0,1).$ 
Тогда   справедливо равенство
$$
((-A)^\beta)^{-\alpha}x=(-A)^{-\beta\alpha}x\ (x\in X_T). \Box 
$$

\end{document}